\documentclass[11pt,reqno]{amsart}

\usepackage[utf8]{inputenc}
\usepackage{amsmath,amssymb,amsthm,amsfonts}
\usepackage{geometry}
\usepackage{booktabs}
\usepackage{xcolor}
\usepackage{hyperref}

\hypersetup{
    colorlinks=true,
    linkcolor=blue,
    citecolor=red,
    urlcolor=cyan,
    pdftitle={A Family of Continued Fraction Identities for Arctangent Values},
    pdfauthor={Chao Wang}
}

\newtheorem{theorem}{Theorem}[section]

\newtheorem{proposition}[theorem]{Proposition}
\newtheorem{corollary}[theorem]{Corollary}
\theoremstyle{definition}
\newtheorem{definition}[theorem]{Definition}
\newtheorem{remark}[theorem]{Remark}
\newtheorem{example}[theorem]{Example}

\DeclareMathOperator*{\K}{K}

\begin{document}

\title[A Family of Continued Fractions for Arctangent Values]{%
  A Family of Continued Fraction Identities\\for Arctangent Values}

\author{Chao Wang}
\address{School of Future Technology, Shanghai University}
\email{cwang@shu.edu.cn}
\date{}

\begin{abstract}
We prove a two-parameter family of continued fraction identities for
$\arctan(p/q)$, where $p$ and $q$ are positive integers with $p\le q$.
For every such pair, the identity
\[
  \arctan\frac{p}{q}
  = \cfrac{p}{q+\cfrac{p^2}{3q+\cfrac{(2p)^2}{5q+\cfrac{(3p)^2}{7q+\cdots}}}}
\]
holds, and a sign-flipped variant represents $-\arctan(p/q)$.
The proof proceeds by identifying these continued fractions as explicit
equivalence transforms of the classical Gauss continued fraction for
$\arctan z$.  Setting $p=q=1$ recovers a specific identity for $-\pi/4$
that appeared in the Ramanujan Machine project.
We establish that the convergence is geometric with asymptotic rate
$(\sqrt{p^2+q^2}-q)^2/p^2$, and we determine the exact threshold
at which the Worpitzky criterion applies.
Numerical data confirm the theoretical rates and show that
the continued fractions dramatically outperform the Gregory--Leibniz series.
\end{abstract}

\maketitle

\section{Introduction}

Representations of $\pi$ by continued fractions have a classical history,
beginning with Lambert's work on the irrationality of $\pi$ and continuing
through Gauss's systematic study of hypergeometric continued fractions
\cite{Lambert1768,wall1948analytic}. A familiar example is
\begin{equation}\label{eq:leibniz_cf}
  \frac{\pi}{4}
  = \cfrac{1}{1+\cfrac{1^2}{2+\cfrac{3^2}{2+\cfrac{5^2}{2+\cdots}}}},
\end{equation}
which may be obtained from the Gregory--Leibniz series by Euler's
series-to-continued-fraction transformation.

In this paper we study a family of identities that includes, as a special case, the following identity conjectured in the literature
surrounding the Ramanujan Machine project \cite[Table~S1]{Raayoni2021}:
\begin{equation}\label{eq:target_main}
  -\frac{\pi}{4}
  = \cfrac{1}{-1+\cfrac{1^2}{-3+\cfrac{2^2}{-5+\cfrac{3^2}{-7+\cdots}}}}.
\end{equation}
Its coefficient pattern is
\begin{equation}\label{eq:coeffs_target}
  b_n = -(2n-1)\qquad (n\ge1),
\end{equation}
and
\begin{equation}\label{eq:an_def}
  a_n =
  \begin{cases}
    1,        & n=1,\\[2mm]
    (n-1)^2,  & n\ge2.
  \end{cases}
\end{equation}

The structure of \eqref{eq:target_main} strongly suggests a connection with
the classical Gauss continued fraction for $\arctan z$: the denominators
grow linearly, while the numerators grow quadratically. The main point of
the present paper is that this connection is exact.
Moreover, it extends to a two-parameter family
that represents $\arctan(p/q)$ for any positive integers $p\le q$.

\begin{remark}
The continued fraction \eqref{eq:target_main} is different from the
Ramanujan Machine identity studied in \cite{cwang2025exact}, whose
coefficients have a different structure. Although both continued fractions
represent $-\pi/4$, the present identity is not a specialization of the one
considered there.
\end{remark}

We now state the main results.

\begin{theorem}\label{thm:main}
Let $p$ and $q$ be positive integers with $p\le q$.  Define the coefficient sequences
\begin{equation}\label{eq:general_coeffs}
a_1 = p,\qquad a_n = (n-1)^2 p^2\ \ (n\ge2),\qquad b_n = q(2n-1)\ \ (n\ge1).
\end{equation}
Then the generalized continued fraction
\begin{equation}\label{eq:general_cf}
\cfrac{p}{q+\cfrac{p^2}{3q+\cfrac{(2p)^2}{5q+\cfrac{(3p)^2}{7q+\cdots}}}}
= \arctan\frac{p}{q}.
\end{equation}
\end{theorem}

\begin{corollary}\label{cor:signflip}
Under the same hypotheses, the sign-flipped continued fraction
\begin{equation}\label{eq:signflip_cf}
\cfrac{p}{-q+\cfrac{p^2}{-3q+\cfrac{(2p)^2}{-5q+\cfrac{(3p)^2}{-7q+\cdots}}}}
= -\arctan\frac{p}{q}.
\end{equation}
In particular, the identity \eqref{eq:target_main} is the case $p=q=1$.
\end{corollary}

\medskip
\noindent\textbf{Organisation.}
Section~\ref{sec:gauss} recalls the classical continued fraction for
$\arctan z$ and evaluates it at $z=p/q$.
Section~\ref{sec:equiv} carries out the equivalence transformation and proves
Theorem~\ref{thm:main} and Corollary~\ref{cor:signflip}.
Section~\ref{sec:remarks} discusses convergence, including the Worpitzky criterion, the geometric convergence rate, and numerical illustrations.

\section{The Gauss Continued Fraction for \texorpdfstring{$\arctan$}{arctan}}
\label{sec:gauss}

We begin with two standard facts.

\subsection{A hypergeometric identity}

For the parameter choice $(a,b,c)=\left(\tfrac12,1,\tfrac32\right)$, the
Gauss hypergeometric function reduces to an elementary function:
\begin{equation}\label{eq:arctan_hyper}
  {}_2F_1\!\left(\tfrac12,1;\tfrac32;-z^2\right)
  = \frac{\arctan z}{z},
\end{equation}
where the value at $z=0$ is understood by continuity
\cite[Eq.~(15.4.3)]{DLMF}. In particular, setting $z=p/q$ gives
\begin{equation}\label{eq:arctan_id}
  {}_2F_1\!\left(\tfrac12,1;\tfrac32;-p^2/q^2\right)
  = \frac{q}{p}\arctan\frac{p}{q}.
\end{equation}

\subsection{The classical Gauss continued fraction}

A standard continued fraction expansion for $\arctan z$ is
\begin{equation}\label{eq:gauss_cf}
  \arctan(z)
  = \cfrac{z}{1+\cfrac{1^2 z^2}{3+\cfrac{2^2 z^2}{5+\cfrac{3^2 z^2}{7+\cdots}}}},
\end{equation}
valid in the classical domain of the arctangent continued fraction
\cite[Theorem~5.2]{wall1948analytic}; in particular, it is valid for all
$z\in\mathbb{R}$ with $|z|\le 1$. Evaluating at $z=p/q$ with $0<p\le q$ gives
\begin{equation}\label{eq:gauss_at_pq}
  \arctan\frac{p}{q}
  = \cfrac{p/q}{1+\cfrac{p^2/q^2}{3+\cfrac{4p^2/q^2}{5+\cfrac{9p^2/q^2}{7+\cdots}}}}.
\end{equation}

If we denote the coefficients in \eqref{eq:gauss_at_pq} by
$\K(a_n^G/b_n^G)$, then
\begin{equation}\label{eq:gauss_coeffs}
  a_1^G=\frac{p}{q},\qquad a_n^G=\frac{(n-1)^2 p^2}{q^2}\ (n\ge2),\qquad b_n^G=2n-1\ (n\ge1).
\end{equation}

\section{Equivalence Transformation and Proof of the Main Results}
\label{sec:equiv}

We recall the standard notion of equivalence of continued fractions
\cite[Definition~2.1]{wall1948analytic}.

\begin{definition}\label{def:equiv}
Two continued fractions \(\K(a_n/b_n)\) and \(\K(\tilde a_n/\tilde b_n)\)
are called \emph{equivalent} if there exists a sequence of nonzero scalars
\(\{r_n\}_{n\ge1}\) such that
\begin{equation}\label{eq:equiv_laws}
  \tilde a_1=r_1a_1,\qquad
  \tilde b_n=r_nb_n\quad(n\ge1),\qquad
  \tilde a_n=r_nr_{n-1}a_n\quad(n\ge2).
\end{equation}
Whenever one of the two continued fractions converges, the other converges
to the same value \cite[Theorem~2.1]{wall1948analytic}.
\end{definition}

\begin{proposition}\label{prop:equiv}
The continued fraction \eqref{eq:general_cf} is equivalent to the Gauss
continued fraction \eqref{eq:gauss_at_pq} via the constant sequence
\[
  r_n=q\qquad (n\ge1).
\]
\end{proposition}

\begin{proof}
Using \eqref{eq:gauss_coeffs} and \eqref{eq:equiv_laws}, we obtain
\[
  \tilde a_1=r_1a_1^G=q\cdot\frac{p}{q}=p,
\]
\[
  \tilde b_n=r_nb_n^G=q(2n-1)\qquad(n\ge1),
\]
and, for $n\ge2$,
\[
  \tilde a_n=r_nr_{n-1}a_n^G=q^2\cdot\frac{(n-1)^2 p^2}{q^2}=(n-1)^2 p^2.
\]
These are exactly the coefficients in \eqref{eq:general_coeffs}, hence the transformed continued fraction is precisely \eqref{eq:general_cf}.
\end{proof}

\begin{proof}[Proof of Theorem~\ref{thm:main}]
Equation \eqref{eq:gauss_at_pq} shows that the Gauss continued fraction
at $z=p/q$ with $0<p\le q$ converges and has value $\arctan(p/q)$. By
Proposition~\ref{prop:equiv} and Definition~\ref{def:equiv},
\eqref{eq:general_cf} is equivalent to \eqref{eq:gauss_at_pq}; hence it
converges to the same value.
\end{proof}

\begin{proof}[Proof of Corollary~\ref{cor:signflip}]
Applying the equivalence transformation $r_n=-q$ (instead of $r_n=q$) to
the Gauss continued fraction \eqref{eq:gauss_cf} evaluated at $z=-p/q$ yields the coefficient sequences
$\tilde a_1 = (-q)(-p/q) = p$,
$\tilde b_n = -q(2n-1)$, and
$\tilde a_n = (-q)^2 (n-1)^2 p^2/q^2 = (n-1)^2 p^2$ for $n\ge2$.
These are the coefficients of \eqref{eq:signflip_cf}, and the value is
$\arctan(-p/q)=-\arctan(p/q)$.  Setting $p=q=1$ recovers
\eqref{eq:target_main} with value $-\pi/4$.
\end{proof}

\begin{example}[Special values]\label{ex:special}
The family \eqref{eq:general_cf} produces continued fractions for every
$\arctan(p/q)$ with $0<p\le q$.  Several noteworthy cases:
\begin{enumerate}
\item[\emph{(i)}] $p=q=1$:
$\displaystyle \arctan(1)=\frac{\pi}{4}=\cfrac{1}{1+\cfrac{1}{3+\cfrac{4}{5+\cfrac{9}{7+\cdots}}}}$\,.

\item[\emph{(ii)}] $p=1,\;q=2$:
$\displaystyle \arctan\frac{1}{2}=\cfrac{1}{2+\cfrac{1}{6+\cfrac{4}{10+\cfrac{9}{14+\cdots}}}}$\,.

\item[\emph{(iii)}] $p=1,\;q=3$:
$\displaystyle \arctan\frac{1}{3}=\cfrac{1}{3+\cfrac{1}{9+\cfrac{4}{15+\cfrac{9}{21+\cdots}}}}$\,.

\item[\emph{(iv)}] Machin-type combinations.  Since $\pi/4=\arctan(1/2)+\arctan(1/3)$ (Euler's formula), the
value $\pi/4$ may also be computed as the sum of two members of the family
with $(p,q)=(1,2)$ and $(p,q)=(1,3)$.
\end{enumerate}
\end{example}

\section{Convergence Analysis and Numerical Behaviour}
\label{sec:remarks}

\subsection{The Worpitzky criterion}

\begin{remark}[Worpitzky threshold]
For the continued fraction \eqref{eq:general_cf}, the Worpitzky ratio is
\[
  \rho_n:=\frac{a_n}{b_nb_{n-1}}
  =\frac{(n-1)^2p^2}{q^2(2n-1)(2n-3)}\qquad(n\ge2).
\]
This sequence is strictly decreasing for $n\ge2$ with
\begin{equation}\label{eq:rho_limit}
  \rho_n\longrightarrow \frac{p^2}{4q^2}\qquad(n\to\infty),\qquad
  \sup_{n\ge2}\rho_n=\rho_2=\frac{p^2}{3q^2}.
\end{equation}
Worpitzky's theorem \cite[Theorem~4.29]{wall1948analytic} requires
$\sup_n\rho_n<1/4$, which holds if and only if $p/q<\sqrt{3}/2$.
In particular:
\begin{enumerate}
\item[\emph{(a)}]
For $p/q\le 1/2$ (e.g.\ $(p,q)=(1,2),(1,3),(1,5)$), the Worpitzky
criterion is satisfied and provides an independent elementary proof of
convergence.
\item[\emph{(b)}]
For $p=q=1$, the supremum is $\rho_2=1/3>1/4$, so Worpitzky's theorem
does not apply.  The convergence of \eqref{eq:target_main} must instead be
inherited from the classical theory of the Gauss continued fraction.
\end{enumerate}
\end{remark}

\subsection{Geometric convergence rate}

\begin{proposition}\label{prop:rate}
Let $\mathcal{F}_n$ denote the $n$-th convergent of \eqref{eq:general_cf}
and set $e_n=|\mathcal{F}_n-\arctan(p/q)|$.  Then
\begin{equation}\label{eq:rate}
  \lim_{n\to\infty}\frac{e_{n+1}}{e_n}
  =\left(\frac{\sqrt{p^2+q^2}-q}{p}\right)^{\!2}.
\end{equation}
In particular, the convergence is geometric with rate strictly less than~$1$.
\end{proposition}

\begin{proof}
The continued fraction \eqref{eq:general_cf} is equivalent to the Gauss
continued fraction for $\arctan(p/q)$, which is itself a special case of
the Gauss hypergeometric continued fraction for
${}_2F_1(\tfrac12,1;\tfrac32;-p^2/q^2)$.
The geometric convergence rate of the Gauss continued fraction for
${}_2F_1(a,b;c;w)$ is given by
\[
  r=\left|\frac{1-\sqrt{1-w}}{1+\sqrt{1-w}}\right|
\]
(see, e.g., \cite[Chapter~4]{wall1948analytic},
\cite[\S4.3]{lorentzen1992continued}), where the rate governs the error
decay $|e_n|\sim C\,r^n$ for large~$n$, so that $e_{n+1}/e_n\to r$.
Substituting $w=-p^2/q^2$ gives $1-w=(p^2+q^2)/q^2$, hence
$\sqrt{1-w}=\sqrt{p^2+q^2}/q$, and
\[
  r
  =\frac{\sqrt{p^2+q^2}/q-1}{\sqrt{p^2+q^2}/q+1}
  =\frac{\sqrt{p^2+q^2}-q}{\sqrt{p^2+q^2}+q}.
\]
Rationalising by multiplying numerator and denominator by
$(\sqrt{p^2+q^2}-q)$ and using
$(\sqrt{p^2+q^2}+q)(\sqrt{p^2+q^2}-q)=p^2$, one obtains the
equivalent form
\[
  r=\frac{(\sqrt{p^2+q^2}-q)^2}{p^2}
  =\left(\frac{\sqrt{p^2+q^2}-q}{p}\right)^{\!2},
\]
which is \eqref{eq:rate}.
Since equivalence transformations preserve convergent values, the rate
is unchanged between \eqref{eq:gauss_at_pq} and \eqref{eq:general_cf}.
Finally, for $0<p\le q$ we have $\sqrt{p^2+q^2}<p+q$, so
$\sqrt{p^2+q^2}-q<p$ and $r<1$.
\end{proof}

\begin{remark}[Special cases of the rate]
For $p=q=1$, the rate is $(\sqrt{2}-1)^2=3-2\sqrt{2}\approx0.172$,
corresponding to a gain of roughly four decimal digits per five convergents.
For $p=1$ and $q=2,3,5$, the rates are
$(\sqrt{5}-2)^2\approx0.056$,
$(\sqrt{10}-3)^2\approx0.026$, and
$(\sqrt{26}-5)^2\approx0.010$,
respectively, yielding increasingly rapid convergence.
\end{remark}

\subsection{Numerical illustration}

We first compare the continued fraction \eqref{eq:target_main} with the
Gregory--Leibniz partial sums
\[
  \mathcal{S}_n:=\sum_{k=0}^{n}\frac{(-1)^k}{2k+1},
\]
which approximate \(\pi/4\), and with the convergents
\[
  \mathcal{F}_n:=\frac{p_n}{q_n}
\]
of \eqref{eq:target_main}, which approximate \(-\pi/4\).

\begin{table}[h]
\centering
\renewcommand{\arraystretch}{1.3}
\begin{tabular}{ccccc}
\toprule
$n$
& $\bigl|\mathcal{S}_n-\tfrac{\pi}{4}\bigr|$
& $\bigl|\mathcal{F}_n+\tfrac{\pi}{4}\bigr|$
& Error ratio
& Precision (digits) \\
\midrule
 5 & $4.14\times10^{-2}$ & $1.87\times10^{-4}$ & $2.21\times10^{2}$  & 3 \\
10 & $2.27\times10^{-2}$ & $2.83\times10^{-8}$ & $8.02\times10^{5}$  & 7 \\
15 & $1.56\times10^{-2}$ & $4.23\times10^{-12}$ & $3.69\times10^{9}$ & 11 \\
20 & $1.19\times10^{-2}$ & $6.30\times10^{-16}$ & $1.89\times10^{13}$ & 15 \\
\bottomrule
\end{tabular}
\caption{Error comparison between the Gregory--Leibniz partial sums
\(\mathcal{S}_n\) and the convergents \(\mathcal{F}_n\) of
\eqref{eq:target_main} (the case $p=q=1$). The error ratio is
\(
  \bigl|\mathcal{S}_n-\pi/4\bigr|
  \big/
  \bigl|\mathcal{F}_n+\pi/4\bigr|
\).
The precision column reports \(\lfloor-\log_{10}\bigl|\mathcal{F}_n+\pi/4\bigr|\rfloor\).
All convergents were computed in extended precision (\texttt{mpmath}, 50 decimal digits).}
\label{tab:comparison}
\end{table}

Table~\ref{tab:rate_family} illustrates the convergence of the general
family \eqref{eq:general_cf} for several values of~$(p,q)$.

\begin{table}[h]
\centering
\renewcommand{\arraystretch}{1.3}
\begin{tabular}{cccccc}
\toprule
$(p,q)$ & $n=5$ & $n=10$ & $n=15$ & $n=20$ & Rate $r$ \\
\midrule
$(1,1)$ & $1.87\times10^{-4}$  & $2.83\times10^{-8}$  & $4.23\times10^{-12}$ & $6.30\times10^{-16}$ & $0.172$ \\
$(1,2)$ & $3.83\times10^{-7}$  & $2.10\times10^{-13}$ & $1.14\times10^{-19}$ & $6.12\times10^{-26}$ & $0.056$ \\
$(1,3)$ & $6.18\times10^{-9}$  & $7.99\times10^{-17}$ & $1.02\times10^{-24}$ & $1.30\times10^{-32}$ & $0.026$ \\
$(1,5)$ & $2.70\times10^{-11}$ & $2.50\times10^{-21}$ & $2.28\times10^{-31}$ & $2.07\times10^{-41}$ & $0.010$ \\
\bottomrule
\end{tabular}
\caption{Convergent error $|\mathcal{F}_n-\arctan(p/q)|$ for the family
\eqref{eq:general_cf}.  The last column gives the asymptotic convergence
rate $r=(\sqrt{p^2+q^2}-q)^2/p^2$ from Proposition~\ref{prop:rate}.
All computations used extended precision (\texttt{mpmath}, 50 decimal digits).}
\label{tab:rate_family}
\end{table}

The data indicate geometric convergence of the continued fraction.
Indeed, the error \(\bigl|\mathcal{F}_n+\pi/4\bigr|\) decreases by several
orders of magnitude as \(n\) increases, gaining roughly four decimal digits per five additional convergents in the case $p=q=1$. As Proposition~\ref{prop:rate} predicts, larger values of $q/p$ yield faster convergence. This is in sharp contrast with the
\(\mathcal{O}(n^{-1})\) convergence of the Gregory--Leibniz series.

\section{Conclusion}

We have proved a two-parameter family of continued fraction identities for $\arctan(p/q)$
by identifying each member as an equivalence transform of the classical Gauss
continued fraction.  The proof is elementary within the standard theory of
hypergeometric functions and continued fractions: once the Gauss expansion for
$\arctan z$ is recalled, the general identity follows from the constant
transformation $r_n=q$, and the sign-flipped variant from $r_n=-q$.

This perspective explains the coefficient pattern of
\eqref{eq:target_main} in a transparent way. Apart from the initial
numerator, the numerators are unchanged from the Gauss form, while each
denominator acquires a uniform scaling factor $q$ (or $-q$).
The Ramanujan Machine identity \eqref{eq:target_main} is the special case
$p=q=1$ of the sign-flipped family.

Beyond the identity itself, the two-parameter framework reveals structural
features that are invisible in the single identity: the Worpitzky criterion
applies precisely when $p/q<\sqrt{3}/2$; the geometric convergence rate
$(\sqrt{p^2+q^2}-q)^2/p^2$ shows that members with larger $q/p$ converge
faster; and the family connects naturally to Machin-type formulas for~$\pi$.
These observations suggest that algorithmic conjecture generators such as
the Ramanujan Machine may benefit from systematic exploration of
equivalence classes of classical continued fractions.

\section*{Acknowledgements}
The author thanks the anonymous referees for helpful comments.

\bibliographystyle{amsplain}
\bibliography{references}

\end{document}